\documentclass{article}%
\usepackage{amssymb}
\usepackage{graphicx}
\usepackage{amsmath}
\usepackage{amsfonts}%
\newtheorem{theorem}{Theorem}

\newtheorem{corollary}[theorem]{Corollary}

\newtheorem{proposition}[theorem]{Proposition}

\newenvironment{proof}[1][Proof]{\textbf{#1.} }{\ \rule{0.5em}{0.5em}}
\begin{document}

\title{The Lie Algebra of the Group of Bisections\\$\vartriangleleft$A Chapter in Synthetic Differential Geometry of
Groupoids$\vartriangleright$}
\author{Hirokazu Nishimura\\Institute of Mathematics, University of Tsukuba\\Tsukuba, Ibaraki, 305-8571, Japan}
\maketitle

\begin{abstract}
Groupoids provide a more appropriate framework for differential geometry than
principal bundles. Synthetic differential geometry is the avant-garde branch
of differential geometry, in which nilpotent infinitesimals are available in
abundance. The principal objective in this paper is to show within our
favorite framework of synthetic differential geometry that the tangent space
of the group of bisections of a microlinear groupoid at its identity is
naturally a Lie algebra. We give essentially distinct two proofs for its
Jacobi identity.

Keywords:generalized gauge transformation, infinitesimal gauge transformation,
topos-theoretic physics, synthetic differential geometry, groupoid, bisection,
Lie algebra, Jacobi identity

\end{abstract}

\section{Introduction and Preliminaries}

\subsection{Groupoids}

Groupoids are, roughly speaking, categories whose morphisms are always
invertible. They provide differential geometry with a more appropriate, more
natural and more general framework than principal bundles. Indeed, principal
bundles can be put down as (a special kind of) groupoids, though not in a
canonical way. Our standard reference on groupoids is MacKenzie \cite{ma1}.
According to Sharpe \cite{Sh1}, Cartan believed that differential geometry was
the study of an Ehresmann connection on a principal bundle. We hold that
differential geometry is the study of an infinitesimal connection on a groupoid.

Given a groupoid $G$ over a base $M$ with its object inclusion map
$\mathrm{id}:M\rightarrow G$ and its source and target projections
$\alpha,\beta:G\rightarrow M$, we denote by $\mathfrak{S}(G)$ the totality of
sections of $\alpha:G\rightarrow M$, i.e., the totality of mappings
$\sigma:M\rightarrow G$ such that $\alpha\circ\sigma$ is the identity mapping
on $M$. Given $\sigma,\rho\in\mathfrak{S}(G)$, $\sigma\ast\rho\in
\mathfrak{S}(G)$ is defined to be
\[
(\sigma\ast\rho)(x)=\sigma((\beta\circ\rho)(x))\rho(x)
\]
for any $x\in M$. It is well known that the set $\mathfrak{S}(G)$ is a monoid
with respect to the operation $\ast$, where the identity is the object
inclusion map $\mathrm{id}:M\rightarrow G$. A section $\sigma$\ of
$\alpha:G\rightarrow M$ is called a \textit{bisection }of $G$\ if $\beta
\circ\sigma$ is a bijection of $M$ onto $M$. We denote by $\mathfrak{B}(G)$
the totality of bisections of $G$. It is well known that the set
$\mathfrak{B}(G)$ is not only a submonoid of $\mathfrak{S}(G)$ but a group by
itself, in which the inverse of $\sigma\in\mathfrak{B}(G)$ is given by
\[
\sigma((\beta\circ\sigma)^{-1}(x))^{-1}%
\]
for any $x\in M$. Providing that $G$ is a Lie groupoid over a
finite-dimensional smooth manifold $M$, the familiar statement that the
totality $\Gamma(\mathcal{A}G)$ of sections of the Lie algebroid
$\mathcal{A}G$ of $G$ is the Lie algebra of the group $\mathfrak{B}(G)$
remains metaphorical at best, unless it is courageous of you to discuss a sort
of infinite-dimensional manifolds seriously or to leave the traditional set theory.

It is well known that if $P$ is a principal bundle over $M$ with its structure
group $H$, and if $G$ is its gauge groupoid $%
\begin{array}
[c]{c}%
\underline{P\times P}\\
H
\end{array}
$, then $\mathfrak{B}(G)$ can be identified with the group of generalized
gauge transformations of $P$.

\subsection{Synthetic Differential Geometry}

We know well that nilpotent infinitesimals were vivid and active in the realm
of burgeoning analysis during the days of Newton and Leibniz and that they
were rampant throughout the works of such pioneers in differential geometry as
Cartan, Lie and Riemann. The 19th century witnessed that so-called sclerotic
$\varepsilon-\delta$ arguments were eviscerating nilpotent infinitesimals
relentlessly and leaving them in oblivion as logicl anathema.

Just as various enlargements of numbers (e.g., from real numbers to complex
numbers) have enriched mathematics tremendously, so does topos theory
originating in sheaf theory and forcing techniques of axiomatic set theory in
the middle of the preceding century. It provides mathematicians with gadgets
for enlarging set theory generously. In particular, it endows differential
geometers with a kind of El Dorado, in which they can do their work with a
cornucopia of nilpotent infinitesimals. Such user-friendly differential
geometry is called \textit{synthetic} differential geometry, for which our
standard reference is Lavendhomme \cite{l1}. To see how to build a
Grothendieck topos for synthetic differential geometry, the reader is referred
to Kock \cite{k1} and Moerdijk and Reyes \cite{mr1}. Our discussion to follow
throughout the rest of the paper is carried out internally within such a
Grothendieck topos. In particular, $G$ and $M$ as well as $\mathbb{R}$
(intended for the set of real numbers with a cornucopia of nilpotent
infinitesimals pursuant to the general Kock-Lawvere axiom) are assumed to be
entities of such a Grothendieck topos and to be microlinear. We denote by $D$
the set $\{d\in\mathbb{R}\mid d^{2}=0\}$ as usual.

\subsection{Vector Fields}

In synthetic differential geometry a vector field can be looked upon in three
distinct but equivalent ways, namely, a section of the tangent bundle, an
infinitesimal flow and an infinitesimal transformation, for which the reader
is referred to pp.69-71 of Lavendhomme \cite{l1}. The equivalence of these
three viewpoints is based on the following elementary exponential laws in any
topos:
\[
(M^{D})^{M}=M^{M\times D}=(M^{M})^{D}%
\]
These three equivalent viewpoints can easily be uplifted to groupoids, as we
will shortly see. The reader should note that since vector fields can not be
identified with differential operators in synthetic differential geometry, it
is by no means trivial to show that vector fields form a Lie algebra.

Given $x\in M$, we denote by $\mathcal{A}_{x}G$ the totality of mappings
$t:D\rightarrow G$ with $t(0)=\mathrm{id}_{x}$ and $(\alpha\circ t)(d)=x$ for
any $d\in D$. We denote by $\mathcal{A}G$ the set-theoretic union of
$\mathcal{A}_{x}G$'s for all $x\in M$. Proceeding set-theoretically, we can
see that $X\in\Gamma(\mathcal{A}G)$ is to be reckoned to be a mapping
$X:M\times D\rightarrow G$ with $X(x,0)=\mathrm{id}_{x}$ and $\alpha
(X(x,d))=x$ for any $x\in M$. Proceeding set-theoretically further, we have
arrived at the most radical viewpoint, which allows us to regard $X\in
\Gamma(\mathcal{A}G)$ as a mapping $X:d\in D\longmapsto X_{d}\in
\mathfrak{S}(G)$ with $X_{0}=\mathrm{id}$, namely, as a tangent vector to
$\mathfrak{S}(G)$ at $\mathrm{id}$. We will see in the next section that the
image of $X$ lies in $\mathfrak{B}(G)\subset\mathfrak{S}(G)$. We will adopt
these equivalent three viewpoints interchangeably. The principal objective in
this paper is to show that $\Gamma(\mathcal{A}G)$ or the tangent space of
$\mathfrak{B}(G)$ at $\mathrm{id}$ forms a Lie algebra with respect to a
geometrically motivated Lie bracket. Kock \cite{k0} has already established,
by exploiting the correspondence between sections of $\mathcal{A}G$ and
vertical right-invariant vector fields on $G$ as orthodox differential
geometers usually do, that sections of $\mathcal{A}G$ form a Lie algebra. We
prefer to do so by uplifting the trinity of the three viewpoints of vector
fields to the context of groupoids, which will give due status to the group
$\mathfrak{B}(G)$. Since we can not identify elements of $\Gamma
(\mathcal{A}G)$ with differential operators in synthetic differential
geometry, it is by no means trivial to establish the Jacobi identity for
$\Gamma(\mathcal{A}G)$. The discussions in this paper degenerate into
classical ones in synthetic differential geometry, providing that $G$ happens
to be $M\times M$ (the pair groupoid of $M$).

\subsection{Strong Differences}

Kock and Lavendhomme \cite{kl} have introduced the notion of strong difference
for microsquares into synthetic differential geometry, a good exposition of
which can be seen in \S 3.4 of Lavendhomme \cite{l1}. Given two microsquares
$\gamma_{+}$ and $\gamma_{-}$ on $M$, their strong difference $\gamma
_{+}\overset{\cdot}{-}\gamma_{-}$ is defined exactly when $\gamma_{+}%
|_{D(2)}=\gamma_{-}|_{D(2)}$, and it is a tangent vector to $M$. Therefore the
operation $\overset{\cdot}{-}$ is a function from $M^{D^{2}}\underset
{M^{D(2)}}{\times}M^{D^{2}}$ to $M^{D}$. The following simple proposition is
the prototype of the general Jacobi identity.

\begin{proposition}
For any $\gamma_{1},\gamma_{2},\gamma_{3}\in M^{D^{2}}$ with $\gamma
_{1}|_{D(2)}=\gamma_{2}|_{D(2)}=\gamma_{3}|_{D(2)}$, we have
\[
(\gamma_{1}\overset{\cdot}{-}\gamma_{2})+(\gamma_{2}\overset{\cdot}{-}%
\gamma_{3})+(\gamma_{3}\overset{\cdot}{-}\gamma_{1})=0
\]

\end{proposition}

In \S 2 of our previous paper \cite{n1} we have relativized the notion of
strong difference $\overset{\cdot}{-}$ for microsquares to get three strong
differences $\underset{1}{\overset{\cdot}{-}}$, $\underset{2}{\overset{\cdot
}{-}}$ and $\underset{3}{\overset{\cdot}{-}}$ for microcubes. More
specifically, the set $M^{D^{3}}$ can be identified with the set
$(M^{D})^{D^{2}}$ via three isomorphisms $\Psi_{i}:M^{D^{3}}\rightarrow
(M^{D})^{D^{2}}\ $($i=$ $1,2,3$) to be defined by
\begin{align*}
\Psi_{1}(\gamma)(d_{1},d_{2})(d)  &  =\gamma(d,d_{1},d_{2})\\
\Psi_{2}(\gamma)(d_{1},d_{2})(d)  &  =\gamma(d_{1},d,d_{2})\\
\Psi_{3}(\gamma)(d_{1},d_{2})(d)  &  =\gamma(d_{1},d_{2},d)
\end{align*}
for any $\gamma\in M^{D^{3}}$ and any $d_{1},d_{2},d\in D$. Let $\Phi
:M^{D\times D}\rightarrow(M^{D})^{D}$ be defined by
\[
\Phi(\gamma)(d_{1})(d_{2})=\gamma(d_{1},d_{2})
\]
for any $\gamma\in M^{D^{2}}$ and any $d_{1},d_{2}\in D$. By way of example,
the strong difference $\overset{\cdot}{-}:(M^{D})^{D^{2}}\underset
{(M^{D})^{D(2)}}{\times}(M^{D})^{D^{2}}\rightarrow(M^{D})^{D}$ for
microsquares on $M^{D}$ can be transferred via isomorphisms $\Psi_{3}$
and$\ \Phi$ to yield the strong difference $\underset{3}{\overset{\cdot}{-}%
}:M^{D^{3}}\underset{D(2)\times D}{\times}M^{D^{3}}\rightarrow M^{D^{2}}$ for
microcubes on $M$. I.e., given microcubes $\gamma_{\pm}$ on $M$, their strong
difference $\gamma_{+}\underset{3}{\overset{\cdot}{-}}\gamma_{-}$ is defined
exactly when $\gamma_{+}|_{D(2)\times D}=\gamma_{-}|_{D(2)\times D}$, and it
is a microsquare on $M$.

We have established the following general Jacobi identity in \cite{n1}. The
general Jacobi identity was revisited, and its proof was elaborated or
reconsidered in \cite{n2} and \cite{n3}.

\begin{theorem}
Let $\gamma_{123},\gamma_{132},\gamma_{213},\gamma_{231},\gamma_{312}%
,\gamma_{321}$ be six microcubes on $M$. As long as the following three
expressions are well defined, they sum up only to vanish:
\begin{align*}
&  (\gamma_{123}\overset{\cdot}{\underset{1}{-}}\gamma_{132})\overset{\cdot
}{-}(\gamma_{231}\overset{\cdot}{\underset{1}{-}}\gamma_{321})\\
&  (\gamma_{231}\overset{\cdot}{\underset{2}{-}}\gamma_{213})\overset{\cdot
}{-}(\gamma_{312}\overset{\cdot}{\underset{2}{-}}\gamma_{132})\\
&  (\gamma_{312}\overset{\cdot}{\underset{3}{-}}\gamma_{321})\overset{\cdot
}{-}(\gamma_{123}\overset{\cdot}{\underset{3}{-}}\gamma_{213})
\end{align*}

\end{theorem}

\section{The Main Result}

\begin{proposition}
\label{t1}Let $X\in\Gamma(\mathcal{A}G)$. For any $(d,d^{\prime})\in D(2)$, we
have
\[
X_{d+d^{\prime}}=X_{d}\ast X_{d^{\prime}}%
\]

\end{proposition}

\begin{proof}
It suffices to note that both sides coincide in case $d=0$ or $d^{\prime}=0$.
\end{proof}

\begin{corollary}
\label{t2}Let $X\in\Gamma(\mathcal{A}G)$. For any $d\in D$, both $\beta\circ
X_{d}$ and $\beta\circ X_{-d}$ are the inverse to each other, so that
$X_{d}\in\mathfrak{B}(G)$.
\end{corollary}

\begin{corollary}
\label{t3}Let $X\in\Gamma(\mathcal{A}G)$. For any $d\in D$, we have
\[
(X_{d})^{-1}=X_{-d}%
\]

\end{corollary}

\begin{proposition}
\label{t4}The group $\mathfrak{B}(G)$ of bisections of $G$ is microlinear. The
$\mathbb{R}$-module $\Gamma(\mathcal{A}G)$ is isomorphic to the tangent space
to $\mathfrak{B}(G)$ at its identity.
\end{proposition}

\begin{proof}
For the former statement, it suffices to note that the set $\mathfrak{B}(G)$
of bisections of $G$ is in bijective correspondence with the set
\[
\{(\varphi,\psi)\in\mathfrak{S}(G)\times\mathfrak{S}(G)\mid\varphi\ast
\psi=\mathrm{id}\text{ and }\psi\ast\varphi=\mathrm{id}\}\text{,}%
\]
which is undoubtedly microlinear. The latter statement follows from Corollary
\ref{t2}.
\end{proof}

\begin{proposition}
\label{t5}Let $X,Y\in\Gamma(\mathcal{A}G)$. The mapping
\[
\lambda:D(2)\rightarrow\mathfrak{B}(G)
\]
which coincides on the axes of $D(2)$ with $X$ and $Y$ respectively is given
by
\[
\lambda(d,d^{\prime})=X_{d}\ast Y_{d^{\prime}}=Y_{d^{\prime}}\ast X_{d}%
\]
Therefore we have
\[
(X+Y)_{d}=X_{d}\ast Y_{d}=Y_{d}\ast X_{d}%
\]

\end{proposition}

\begin{proof}
It suffices to note that both mappings $\lambda^{\prime},\lambda^{\prime
\prime}:D(2)\rightarrow\mathfrak{B}(G)$ given by
\begin{align*}
\lambda^{\prime}(d,d^{\prime})  &  =X_{d}\ast Y_{d^{\prime}}\\
\lambda^{\prime\prime}(d,d^{\prime})  &  =Y_{d^{\prime}}\ast X_{d}%
\end{align*}
coincide on the axes of $D(2)$ with $X$ and $Y$ respectively.
\end{proof}

Given $X,Y\in\Gamma(\mathcal{A}G)$, we define $\lambda:D^{2}\rightarrow
\mathfrak{B}(G)$ to be
\[
\lambda(d_{1},d_{2})=Y_{-d_{2}}\ast X_{-d_{1}}\ast Y_{d_{2}}\ast X_{d_{1}}%
\]
for any $(d_{1},d_{2})\in D(2)$. It is easy to see that
\begin{align*}
\lambda(d,0)  &  =Y_{0}\ast X_{-d_{1}}\ast Y_{0}\ast X_{d_{1}}=X_{-d_{1}}\ast
X_{d_{1}}=\mathrm{id}\\
\lambda(0,d)  &  =Y_{-d_{2}}\ast X_{0}\ast Y_{d_{2}}\ast X_{0}=Y_{-d_{2}}\ast
Y_{d_{2}}=\mathrm{id}%
\end{align*}
for any $d\in D$. Therefore, by Proposition 3 in \S 2.3 of \cite{l1}, there
exists a unique element of $\Gamma(\mathcal{A}G)$, called the \textit{Lie
bracket} of $X$ and $Y$ and denoted by $[X,Y]$, such that for any $d_{1}%
,d_{2}\in D$, $[X,Y]_{d_{1}d_{2}}\ $is
\[
Y_{-d_{2}}\ast X_{-d_{1}}\ast Y_{d_{2}}\ast X_{d_{1}}%
\]

Now we come to the main result of this paper.

\begin{theorem}
The $\mathbb{R}$-module $\Gamma(\mathcal{A}G)$ is a Lie algebra with respect
to $[\cdot,\cdot]$.
\end{theorem}

\begin{proof}
It suffices to establish that for any $a\in\mathbb{R}$ and any $X,Y,Z\in
\Gamma(\mathcal{A}G)$, we have
\begin{align}
\lbrack aX,Y]  &  =a[X,Y]\label{1.1}\\
\lbrack X+Y,Z]  &  =[X,Z]+[Y,Z]\label{1.2}\\
\lbrack X,Y]  &  =-[Y,X]\label{1.3}\\
\lbrack X,[Y,Z]]+[Y,[Z,X]]+[Z,[X,Y]]  &  =0 \label{1.4}%
\end{align}
Here we deal only with the first three identities, leaving the last \ identity
to the succeeding two sections. For (\ref{1.1}) we have to note that
\begin{align*}
&  [aX,Y]_{_{d_{1}d_{2}}}\\
&  =Y_{-d_{2}}\ast(aX)_{-d_{1}}\ast Y_{d_{2}}\ast(aX)_{d_{1}}\\
&  =Y_{-d_{2}}\ast X_{-ad_{1}}\ast Y_{d_{2}}\ast X_{ad_{1}}\\
&  =[X,Y]_{a_{d_{1}d_{2}}}%
\end{align*}
The identity (\ref{1.2}) follows directly from (\ref{1.1}) by dint of
Proposition 10 of \S 1.2 of Lavendhomme \cite{l1}. For (\ref{1.3}) It suffices
to show equivalently that
\[
\lbrack X,Y]+[Y,X]=0\text{,}%
\]
which follows from the following computation:
\begin{align*}
&  ([X,Y]+[Y,X])_{d_{1}d_{2}}\\
&  =[X,Y]_{_{d_{1}d_{2}}}\ast\lbrack Y,X]_{_{d_{1}d_{2}}}\text{ \ \ (by
Proposition \ref{t5})}\\
&  =Y_{-d_{2}}\ast X_{-d_{1}}\ast Y_{d_{2}}\ast X_{d_{1}}\ast X_{-d_{1}}\ast
Y_{-d_{2}}\ast X_{d_{1}}\ast Y_{d_{2}}\\
&  =\mathrm{id}%
\end{align*}

\end{proof}

\section{The First Approach to the Jacobi Identity}

We will mimic Nishimura's \cite{n0} shortest and simplest proof of the Jacobi
identity of vector fields in synthetic differential geometry. The basic idea
therein is that the familiar relationship between the Lie bracket of two
vector fields and the Lie derivative of a vector field with respect to the
other reduces the well-known Jacobi identity of vector fields to something
akin to the Leibniz rule in elementary differential calculus.

Given $X\in\Gamma(\mathcal{A}G)$ and $\sigma\in\mathfrak{B}(G)$, we denote by
$\sigma_{\ast}X$ the tangent vector $d\in D\longmapsto$ $\sigma\ast X_{d}%
\ast\sigma^{-1}\in\mathfrak{B}(G)$ to $\mathfrak{B}(G)$ at $\mathrm{id}$.
Given $X,Y\in\Gamma(\mathcal{A}G)$, we define the \textit{Lie derivative of}
$Y$ \textit{with respect to} $X$, denoted by $\mathbf{L}_{X}Y$, to be the
unique tangent vector to $\mathfrak{B}(G)$ at $\mathrm{id}$ such that
\[
(X_{-d})_{\ast}Y-Y=d\mathbf{L}_{X}Y
\]
for any $d\in D$. The relationship between $\mathbf{L}_{X}Y$ and $[X,Y]$ is as
simple as follows:

\begin{theorem}
\label{t2.1}For any $X,Y\in\Gamma(\mathcal{A}G)$, we have $\mathbf{L}%
_{X}Y=[X,Y]$.
\end{theorem}

\begin{proof}
For any $d,d^{\prime}\in D$, we have
\begin{align*}
&  ((X_{-d})_{\ast}Y-Y)_{d^{\prime}}\\
&  =X_{-d}\ast Y_{d^{\prime}}\ast X_{d}-Y_{d^{\prime}}\\
&  =Y_{-d^{\prime}}\ast X_{-d}\ast Y_{d^{\prime}}\ast X_{d}\\
&  \text{[by Proposition \ref{t5}]}\\
&  =[X,Y]_{dd^{\prime}}\\
&  =(d[X,Y])_{d^{\prime}}%
\end{align*}
Therefore
\[
(X_{-d})_{\ast}Y-Y=d[X,Y]
\]
for any $d\in D$, which is tantamount to saying that $\mathbf{L}_{X}Y=[X,Y]$.
\end{proof}

With the above theorem in mind, we can say that the following theorem is none
other than Jacobi's identity of $\Gamma(\mathcal{A}G)$ in disguise.

\begin{theorem}
\label{t2.2}For all $X,Y,Z\in\Gamma(\mathcal{A}G)$, we have
\begin{equation}
\mathbf{L}_{X}[Y,Z]=[\mathbf{L}_{X}Y,Z]+[Y,\mathbf{L}_{X}Z] \label{2.2}%
\end{equation}

\end{theorem}

\begin{proof}
It is easy to see that for any $d\in D$, we have
\[
(X_{-d})_{\ast}[Y,Z]=[(X_{-d})_{\ast}Y,(X_{-d})_{\ast}Z]
\]
Therefore we have
\begin{align*}
&  (X_{-d})_{\ast}[Y,Z]-[Y,Z]\\
&  =[(X_{-d})_{\ast}Y,(X_{-d})_{\ast}Z]-[Y,Z]\\
&  =[d\mathbf{L}_{X}Y+Y,d\mathbf{L}_{X}Z+Z]-[Y,Z]\\
&  =d([\mathbf{L}_{X}Y,Z]+[Y,\mathbf{L}_{X}Z])
\end{align*}
Consequently we have
\[
\mathbf{L}_{X}[Y,Z]=[\mathbf{L}_{X}Y,Z]+[Y,\mathbf{L}_{X}Z]\text{,}%
\]
as was expected.
\end{proof}

The reader should note that the proof of the above theorem is no more
difficult than the synthetic proof of Leibniz's rule in elementary
differential calculus, which is to be seen e.g., in Proposition 1 of \S 1.2 of
Lavendhomme \cite{l1}. Now Jacobi's identity in $\Gamma(\mathcal{A}G)$ with
respect to $[\cdot,\cdot]$ is an easy consequence of the preceding theorem, as
we will see just below.

\begin{proof}
(The Jacobi identity) It suffices to establish the following version of the
Jacobi identity:
\[
\lbrack X,[Y,Z]]=[[X,Y],Z]+[Y,[X,Z]]
\]
for all $X,Y,Z\in\Gamma(\mathcal{A}G)$. Since we have
\begin{align*}
\lbrack X,[Y,Z]]  &  =\mathbf{L}_{X}[Y,Z]\\
\lbrack X,Y]  &  =\mathbf{L}_{X}Y\\
\lbrack X,Z]  &  =\mathbf{L}_{X}Z
\end{align*}
by dint of Theorem \ref{t2.1}, the desired equality is no other than a
reformulation of (\ref{2.2}) of Theorem \ref{t2.2}.
\end{proof}

\section{The Second Approach to the Jacobi Identity}

Given $X_{1},...,X_{n}\in\Gamma(\mathcal{A}G)$, we define $n$-mircrocube
$X_{n}\circledast...\circledast X_{1}$ on $\mathfrak{B}(G)$ to be
\[
(X_{n}\circledast...\circledast X_{1})(d_{1},...d_{n})=(X_{n})_{d_{n}}%
\ast...\ast(X_{1})_{d_{1}}%
\]
for any $(d_{1},...d_{n})\in D^{n}$. Given an $n$-mircrocube $\gamma$ on
$\mathfrak{B}(G)$ and a permutation $\varepsilon$ of the first $n$ natural
numbers, we define another $n$-mircrocube $\sum_{\varepsilon}\gamma$ on
$\mathfrak{B}(G)$ to be
\[
(\sum_{\varepsilon}\gamma)(d_{1},...d_{n})=\gamma(\varepsilon(d_{1}%
),...\varepsilon(d_{n}))
\]
for any $(d_{1},...d_{n})\in D^{n}$. A cycle $\varepsilon$ of length $k$ is
usually denoted by $(j,\varepsilon(j),\varepsilon^{2}(j),...,\varepsilon
^{k-1}(j))$, where $j$ is not fixed by $\varepsilon$.

\begin{proposition}
Let $X,Y\in\Gamma(\mathcal{A}G)$. Both $Y\circledast X$ and $%
{\textstyle\sum}
(X\circledast Y)$ coincide on the axes, and we have
\[
\lbrack X,Y]=Y\circledast X\overset{\cdot}{-}%
{\textstyle\sum_{(12)}}
(X\circledast Y)
\]

\end{proposition}

\begin{proof}
Since we have
\begin{align*}
(Y\circledast X)(d_{1},d_{2})  &  =Y_{d_{2}}\ast X_{d_{1}}\\%
{\textstyle\sum_{(12)}}
(X\circledast Y)(d_{1},d_{2})  &  =X_{d_{1}}\ast Y_{d_{2}}%
\end{align*}
for any $(d_{1},d_{2})\in D^{2}$, it is easy to see that
\begin{align*}
(Y\circledast X)(d,0)  &  =X_{d}=%
{\textstyle\sum_{(12)}}
(X\circledast Y)(d,0)\\
(Y\circledast X)(0,d)  &  =Y_{d}=%
{\textstyle\sum_{(12)}}
(X\circledast Y)(0,d)
\end{align*}
for any $d\in D$, so that the right-hand side of the desired equality is
meaningful. We define a function $\lambda:D^{3}\{(1,3),(2,3)\}\rightarrow
\mathfrak{B}(G)$ to be
\[
\lambda(d_{1},d_{2},d_{3})=X_{d_{1}}\ast\lbrack X,Y]_{d_{3}}\ast Y_{d_{2}}%
\]
Then it is easy to see that
\begin{align*}
\lambda(d_{1},d_{2},0)  &  =X_{d_{1}}\ast Y_{d_{2}}=%
{\textstyle\sum_{(12)}}
(X\circledast Y)(d_{1},d_{2})\\
\lambda(d_{1},d_{2},d_{1}d_{2})  &  =X_{d_{1}}\ast\lbrack X,Y]_{d_{1}d_{2}%
}\ast Y_{d_{2}}\\
&  =X_{d_{1}}\ast\lbrack Y,X]_{d_{1}(-d_{2})}\ast Y_{d_{2}}\\
&  =X_{d_{1}}\ast(X_{-d_{1}}\ast Y_{d_{2}}\ast X_{d_{1}}\ast Y_{-d_{2}})\ast
Y_{d_{2}}\\
&  =Y_{d_{2}}\ast X_{d_{1}}\\
&  =(Y\circledast X)(d_{1},d_{2})
\end{align*}
Therefore we have
\begin{align*}
&  (Y\circledast X\overset{\cdot}{-}%
{\textstyle\sum_{(12)}}
(X\circledast Y))(d)\\
&  =\lambda(0,0,d)\\
&  =[X,Y]_{d}%
\end{align*}
for any $d\in D$, which completes the proof.
\end{proof}

\begin{proposition}
Given $X,Y,Z\in\Gamma(\mathcal{A}G)$, we define six microcubes on
$\mathfrak{B}(G)$ as follows:
\begin{align*}
\gamma_{123}  &  =Z\circledast Y\circledast X\\
\gamma_{132}  &  =%
{\textstyle\sum_{(23)}}
(Y\circledast Z\circledast X)\\
\gamma_{213}  &  =%
{\textstyle\sum_{(12)}}
(Z\circledast X\circledast Y)\\
\gamma_{231}  &  =%
{\textstyle\sum_{(123)}}
(X\circledast Z\circledast Y)\\
\gamma_{312}  &  =%
{\textstyle\sum_{(132)}}
(Y\circledast X\circledast Z)\\
\gamma_{321}  &  =%
{\textstyle\sum_{(13)}}
(X\circledast Y\circledast Z)
\end{align*}
Then the right-hand sides of the following three identities are meaningful,
and all the three identities hold:
\begin{align*}
\lbrack X,[Y,Z]]  &  =(\gamma_{123}\overset{\cdot}{\underset{1}{-}}%
\gamma_{132})\overset{\cdot}{-}(\gamma_{231}\overset{\cdot}{\underset{1}{-}%
}\gamma_{321})\\
\lbrack Y,[Z,X]]  &  =(\gamma_{231}\overset{\cdot}{\underset{2}{-}}%
\gamma_{213})\overset{\cdot}{-}(\gamma_{312}\overset{\cdot}{\underset{2}{-}%
}\gamma_{132})\\
\lbrack Z,[X,Y]]  &  =(\gamma_{312}\overset{\cdot}{\underset{3}{-}}%
\gamma_{321})\overset{\cdot}{-}(\gamma_{123}\overset{\cdot}{\underset{3}{-}%
}\gamma_{213})
\end{align*}

\end{proposition}

Now it is easy to see that the Jacobi identity in $\Gamma(\mathcal{A}G)$
follows readily from the general Jacobi identiy.

\end{document}